\documentclass[11pt]{amsart}

\usepackage{graphicx}
\usepackage{url}

\title{Leonardo and the Pseudo-RCO}

\author[J.M.~Wills]{J\"org M.~Wills}
\address{University of Siegen, Fachbereich Mathematik, Emmy-Noether-Campus, Walter-Flex-Str. 3, 57068 Siegen}
\email{wills@mathematik.uni-siegen.de}

\begin{document}

\begin{abstract}
There are two hypotheses on Leonardo's polyhedron based on the
Pseudo-RCO and  drawn for Luca Pacioli's book: Leonardo made an
error, or: Leonardo draw it with intention, as it is.
We give arguments, which support the Intention-hypothesis. 
\end{abstract}

\maketitle

\sloppy

In 2011 the Dutch artist and mathematician Rinus Roelofs discovered
an  ``error'' in one of the famous polyhedra paintings for Luca Pacioli's
book ``De Divina proportione'' \cite{pacioli}.
Roelofs' discovery appeared in  the Scientific American \cite{huylebrouck} and various
other scientific journals.
The polyhedron is rather complicate (see the figure; the best figures are in
Carlo H.Sequin's contribution~\cite{sequin}).

\medskip

\begin{center}
\includegraphics[width=7cm]{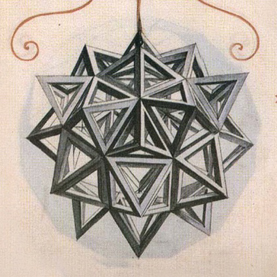}

{\small Leonardo's polyhedron based on the Pseudo-RCO}
\end{center}

\medskip

The ``error'' is not easy to find, which explains
its discovery after more than 500~years.
But the crucial point is, that the polyhedron is correctly drawn:
It does not contain any false line or vertex. It  is only an unexpected
polyhedron, based on the Pseudo-RCO rather than the RCO (rhombi-cuboctahedron), 
one of the 13~Archimedean solids.
Leonardo's job was to draw the 5~Platonic solids and 13~Archimedean solids
for Pacioli's book and for each  one four variations: The basic solid,
the  corresponding star polyhedron (pyramids on the faces) and for
both the  edge-skeleton, where the edges are replaced by thin struts.
In each case the basic solid is the simplest and the star polyhedron with struts
the most complicated polyhedron. For shortness we call it the final polyhedron.
If one looks at the complicated final polyhedron (see the figure), one might
get the impression, that Leonardo lost the track. But:
The construction of the final polyhedron had to start with the basic
polyhedron, and then it was built up step by step. During this process
the basis could not change (like the foundation of a house).
So Leonardo made his decision (or error) at the beginning with the
basic polyhedron. It is rather difficult to mix up the RCO and the
Pseudo RCO at this step: If he started with a model (of 18~squares and
8~regular triangles) or if he started with a drawing (of the edge-skeleton
of a cube and successive edge-cuts and vertex-cuts).
In both cases it is nearly impossible to  ignore the  symmetric RCO
and to choose the  less symmetric Pseudo-RCO.
This supports the ``Intention-Hypothesis'', or as Carlo H.Sequin formulated~\cite{sequin}:
``Leonardo knew, what he was doing''.
Finally one may ask for Leonardo's motivation: Why did he draw the
unexpected polyhedron? 
Here we leave science and come to psychology and speculation:
Perhaps he wanted to irritate the viewer. Or he wanted to break the routine,
because he had to draw roughly 60~polyhedra, perhaps too much routine
for a genius. We do not know the motivation.

\end{document}